\documentclass[11pt]{article}
\usepackage{a4}
\usepackage{graphicx}
\usepackage{parskip}

\usepackage{amsfonts}
\usepackage{natbib}
\usepackage{algorithm}
\usepackage{algorithmic}

\usepackage{latexsym}

\newtheorem{theorem}{Theorem}[section]

\newtheorem{lemma}{Lemma}[section]

\newtheorem{remark}{Remark}[section]
\newtheorem{example}{Example}[section]
\newtheorem{definition}{Definition}[section]

\newcommand{\R}{\mathbb{R}}

\newcommand{\argmin}{\mathop{\arg \min}\limits}

\newcommand{\EE} {{\rm I\hskip-0.48em E}}

\usepackage{natbib}

\begin{document}

\begin{center}
\LARGE{Weakly decomposable regularization penalties and structured sparsity}\\\vspace{1.5cm}
 \normalsize{\textit{Weakly decomposable regularization penalties}}
\vspace{1.5cm}

\Large{Sara van de Geer \\Seminar
  f\"ur Statistik, ETH Zurich}
\end{center}

\newpage

{\bf Abstract} It has been shown in literature that the Lasso estimator,
or $\ell_1$-penalized least squares estimator, enjoys good oracle properties.
This paper examines which special properties of
the $\ell_1$-penalty allow for sharp oracle results, and
then extends the situation to general norm-based penalties
that satisfy a weak decomposability condition.

\noindent\textit{Key words}: {Lasso, sharp oracle inequality, weakly decomposable norm, sparsity, regularization}

\vspace*{.3in}

\section{Introduction} \label{introduction.section}

The Lasso (\cite{tibs96}) has become extremely popular in the last several years. It is a computationally tractable method for
high-dimensional models, with good theoretical properties. 
Several types of modifications of the Lasso have been introduced and studied, such as
the fused Lasso (\cite{Tibshirani:05}) and  the smoothed Lasso (\cite{hebirivdg11}).
In this paper, we are primarily interested in extensions of the $\ell_1$-penalty to general structured sparsity penalties
such as the group Lasso
introduced by \cite{lin06grouplasso}) and further structured versions given by \cite{Zhao:2009}, \cite{Jacob:2009},
\cite{Jenatton:2009} and \cite{Micchelli:2010}. We will  provide sharp versions 
of the oracle inequalities given in \cite{Bach:10} 
and extend the sharp oracle results 
for the Lasso and nuclear norm penalization as given in
\cite{koltchinskii2011nuclear} and \cite{koltchinskii2011oracle} to general structured sparsity penalties, where we
in addition prove inequalities for the estimation error.

Consider the linear model
$$Y = X \beta^0 + \epsilon , $$
where $Y$ is an $n$-vector of observations, $X$ is a given $n \times p$ matrix, $\epsilon$ is an $n$-vector of errors and
$\beta^0 $ is a $p$-vector of unknown coefficients.
The Lasso estimator is
$$\hat \beta := \arg \min_{\beta \in \R^p} \biggl \{ \| Y - X \beta \|_n^2 + 2 \lambda \| \beta \|_1 \biggr \} . $$
Here, $\| \beta \|_1:= \sum_{j=1}^p | \beta_j | $ denotes the $\ell_1$-norm of the vector $\beta$ and 
for a vector $v \in \R^n$ we let $\| v \|_n$ be the normalized Euclidean norm
$ \| v \|_n := \sqrt { v^T v/n } $. 
Finally $\lambda > 0 $ is a tuning parameter. The $\ell_1$-penalty is a variable selection
or (soft-)thresholding type penalty: the larger $\lambda$, the more coefficients
$\hat \beta_j$ will be set to zero.

In this paper, we first briefly review  a sharp oracle result  of \cite{koltchinskii2011nuclear}
for the Lasso estimator.
We then extend the sharp oracle result to other norm-penalties, satisfying
a weak decomposability condition as given in Section \ref{separability.section}.

The paper is organized as follows. We first introduce the concept ``effective sparsity" 
in Section \ref{compatibility.section}. Effective sparsity plays a crucial role
in all our results. As a benchmark, we then  restate in Section \ref{Lasso.section}
an oracle inequality from \cite{koltchinskii2011nuclear}
for the Lasso. Theorem \ref{main.theorem}  in Section \ref{separability.section} 
contains the main result. It extends the
$\ell_1$-norm penalty to general weakly decomposable norm-penalties. Some examples are given
in Section \ref{examples.section}. In Section \ref{eigenvalues.section} we consider comparison of
the effective sparsity based on  the $\ell_1$-norm to the effective sparsity
based on a different norm. A brief discussion of the results and further research is
given in Section \ref{discussion}. Finally, Section \ref{proofs.section} contains the proofs.

\section{The $\ell_1$-eigenvalue and effective sparsity for the $\ell_1$-norm}
\label{compatibility.section}
To state an oracle result, we need to define the $\ell_1$-eigenvalue 
$\delta(L,S) $, where
$L > 0$ is a constant and 
$S \subset \{ 1 , \ldots , p \}$ is an index set. 
We use the notation
$$\beta_{j,S} := \beta_j {\rm l} \{ j \in S \} , \ j=1 , \ldots , p. $$
Thus $\beta_S$ is a $p$-vector with zero entries at the indexes $j \notin S$.
We will sometimes identify $\beta_S$ with the vector
$\{ \beta_j \}_{j \in S} \in \R^{|S|}$.

\begin{definition} For constant $L>0$ and an index set $S$, the $\ell_1$-eigenvalue  is
$$ \delta (L,S) := \min \biggl \{  \| X \beta_S - X \beta_{S^c} \|_n : \ \| \beta_S \|_1 = 1 , \ \| \beta_{S^c } \|_1 \le L \biggr \}. $$
The
compatibility
constant is
$$\phi^2(L,S) := |S| \delta^2 (L,S) . $$
\end{definition} 

 The geometric interpretation 
of the $\ell_1$-eigenvalue, 
as given in \cite{Lederer:11}, is as follows. Let $X_j \in \R^n $ denote the $j$-th column of $X$
($j=1 , \ldots , p $). 
The set $\{ X \beta_S : \ \| \beta_S \|_1 =1 \} $
is the convex hull of the vectors $\{ \pm X_j \}_{j \in S} $ in $\R^n$. Likewise, the set
$\{ X \beta_{S^c} : \ \| \beta_{S^c} \|_1 \le L \} $ is the convex hull including interior
of the vectors $\{ \pm L X_j \}_{j \in S^c} $. Thus, the $\ell_1$-eigenvalue 
$\delta (L, S)$
is the distance
between these two sets.  
We note that:\\
- if $L$ is large the $\ell_1$-eigenvalue will be small,\\
- it will also be small if  the vectors in $S$ exhibit strong correlation with those in $S^c$,\\
- when the vectors in $\{ X_j \}_{j \in S}$ are linearly dependent, it holds
that 
$$\{ X \beta_{S} : \ \| \beta_S \|_1 =1 \} = \{ X \beta_{S} : \ \| \beta_S \|_1 \le 1 \},$$ 
and hence then
$\delta (L,S) =0 $. \\
The compatibility constant was introduced in \cite{vandeG07}. Its name comes from the idea that 
when $\phi (L,S)$ is large the
normalized Euclidean norm $\| \cdot \|_n$ and the $\ell_1$-norm $\| \cdot \|_1$ are in sense
compatible. The difference between the compatibility constant and the squared
$\ell_1$-eigenvalue lies only in the normalization by the size $|S|$ of the set $S$. This normalization is inspired by
the orthogonal case, which we detail in the following example.

\begin{example} Suppose that the columns of $X$ are all orthogonal:
$X_j^T X_k=0$ for all $j \not=k $. Assume moreover the normalization
$\| X_j \|_n =1$ for all $j$. Then clearly,
$$\| X\beta_S - X \beta_{S^c} \|_n = \| \beta_S - \beta_{S^c} \|_2 , $$
where $\| \beta \|_2 := \sqrt {\sum_{j=1}^p \beta_j^2 }$ is the
$\ell_2$-norm of the vector $\beta$. But
$$\| \beta_S - \beta_{S^c} \|_2^2 = \| \beta_S \|_2^2 + \| \beta_{S^c} \|_2^2 \ge
\| \beta_S \|_2^2 \ge \| \beta_S \|_1^2 /|S|,$$
and in fact
$$\min_{\| \beta_{S^c} \|_1 \le L, \ \| \beta_S \|_1 = 1} \| \beta_S - \beta_{S^c} \|_2^2 = 
\min_{\| \beta_S \|_1 = 1 } \| \beta_S \|_2^2 = 1/|S| . $$
It follows that $\delta^2 (L,S) =1/ |S| $ and $\phi^2(L,S) = 1 $.
\end{example}

A vector $\beta$ is called sparse if it has only few non-zero coefficients.
That is, the cardinality $|S_{\beta} |$ of the set $S_{\beta} := \{ j:\ \beta_j \not= 0 \}$ is small.
We call $|S_{\beta}|$ the sparsity-index of $\beta$. More generally, we call $|S|$ the sparsity index of the set $S$.
The {\it effective} sparsity, as defined in \cite{Mueller:11}, takes into account the correlation structure in the design matrix
$X$. 

\begin{definition} For a set $S$ and constant $L>0$, the effective sparsity $\Gamma^2(L,S)$ is the inverse
of the squared $\ell_1$-eigenvalue, that is
$$\Gamma^2 (L,S) = { 1 \over \delta^2 (L,S) } . $$
\end{definition}

In other words, for orthogonal design the effective sparsity of a set $S$ is its cardinality,
and in general, it is the inverse of the squared distance between the convex hull
$\{ X \beta_S: \ \| \beta_S \|_1 = 1 \} $ and the convex set $\{ X \beta_{S^c} : \ 
\| \beta_{S^c } \|_1 \le L \} $. 

Finally, we give a small numerical example from \cite{Mueller:11}.

\begin{example} As a simple numerical example, let us suppose $n=2$, $p=3$, $S= \{ 3 \}$, 
and
$${X} = \sqrt {n}\pmatrix { 5/13 & 0 & 1 \cr 12/13 & 1 & 0 \cr}   . $$
Since the sparsity index is $|S|=1$, the $\ell_1$-eigenvalue $\delta(L,S)$ is equal to 
the square root $\phi(L,S)$ of the compatibility constant, and 
equal to the distance of ${ X}_1$ to line that connects $L {X}_1 $ and $-L {X}_2 $, that is
$$ \delta (L,S) = \max \{ (5-L)/\sqrt {26} , 0 \} . $$
Hence, for example for $L=3$ the effective sparsity is $\Gamma^2 (3, S) = 13/2$.\\
Alternatively, when
$${ X} =\sqrt {n}  \pmatrix {12/13 & 0 & 1 \cr 5/13 & 1 & 0 \cr  } ,$$
then for example $\delta(3,S)=0$ and hence $\Gamma^2 (3,S) = \infty$. This is due to the sharper angle between ${X}_1$ and ${ X}_3$.

\end{example}

\section{An oracle inequality for the $\ell_1$-norm}\label{Lasso.section}

For a vector
$w \in \R^p$, we let
$\| w \|_{\infty} := \max_{1 \le j \le p } |w_j| $ be the uniform norm.
The following theorem is a slight extension of \cite{koltchinskii2011nuclear} (we use the effective
sparsity instead of restricted eigenvalues).
The  sparsity oracle inequality in this theorem is a simple consequence of the following properties of the $\ell_1$-norm:\\
$\bullet$ { Dual norm equality:} $ \sup \{ | w^T \beta |  : \ \| \beta \|_1 \le 1 \} =\| w \|_{\infty}$, $\forall \ w$,\\
$\bullet$ { Triangle inequality :} $\| \beta  + \tilde \beta \|_1 \le  \|  \beta \|_1 + \|  \tilde \beta \|_1 $,  \ $\forall \ \beta , \tilde \beta $,\\
$\bullet$ { Decomposability:} $\| \beta \|_1 = \| \beta_{S} \|_1 + \| \beta_{S^c} \|_1$,  $\forall \ \beta, \ S$. \\

Note that the triangle inequality implies convexity: $\| \alpha \beta + (1- \alpha) \tilde \beta \|_1 \le
\alpha \| \beta \|_1 + (1- \alpha) \| \tilde \beta \|_1 $, $\forall \ \beta \tilde \beta $ and all $0 \le \alpha \le 1 $.
Convexity of the penalty is crucial for deriving oracle inequalities that are sharp. Lemma
\ref{convexpenalty.lemma} gives the details. 

Recall the notation
$$S_{\beta} := \{ j: \ \beta_j \not= 0 \} , \ \beta \in \R^p . $$

\begin{theorem}\label{leastsquares} (\cite{koltchinskii2011nuclear})
Let for $S \subset \{ 1 , \ldots , p \}$
$$\lambda^S := \| (\epsilon^T X )_S \|_{\infty} /n , \ \lambda^{S^c} := \| (\epsilon^T X )_{S^c} \|_{\infty} /n .  $$ 
Define  for $\lambda >  \lambda^{S^c}  $
$$L_S:= { \lambda + \lambda^S \over \lambda- \lambda^{S^c} } . $$
Then
$$ \| X ( \hat \beta - \beta^0 ) \|_n^2  \le \min_{\beta \in \R^p,\ S= S_{\beta} , \ \lambda > \lambda^{S^c}    } \biggl \{ \| X( \beta - \beta^0) \|_n^2 +  (\lambda + \lambda^S )^2 \Gamma^2 (L_S, S  ) \biggr \} .  $$
\end{theorem}

Thus, the Lasso trades off an approximation error $\| X ( \beta - \beta^0) \|_n^2 $ with
an estimation error $(\lambda + \lambda^S )^2 \Gamma^2 (L, S_{\beta}  ) $. 
The above oracle inequality is called sharp because the constant in front of the approximation error
$\| X ( \beta - \beta^0 )\|_n^2 $ is one. Apart from \cite{koltchinskii2011nuclear} and \cite{koltchinskii2011oracle}, 
results in literature are mostly non-sharp versions, 
with a constant larger than one
in front of the approximation error, see e.g.\ \cite{BvdG2011}.  It is interesting to note that
convexity of the penalty plays a crucial role, e.g., with the $\ell_0$-penalty one
cannot arrive at sharp oracle results. 
Observe that we do not present a bound for the $\ell_1$-error in Theorem \ref{leastsquares}.
We will show how such a bound can be included in the results 
in Theorem \ref{main.theorem}.

\begin{remark} \label{repeat.remark}
It is as yet not clear to what extent $\ell_1$-eigenvalue
conditions are necessary for oracle behavior of the prediction error
$\| X ( \hat \beta - \beta^0 ) \|_n^2 $ of the Lasso estimator.
For example, if the design matrix $X$ has repeated columns
(or columns that are proportional) in the set $S$, then the
$\ell_1$-eigenvalue will be zero. A reparametrization argument shows however that
the Lasso estimator behaves as if repeated columns are treated as one.
\end{remark}

\section{A sharp oracle inequality for general weakly decomposable penalties}\label{separability.section}

Let $\Omega$ be some norm on $\R^p$, and let $\hat \beta$ be the norm-penalized estimator
$$\hat \beta := \hat \beta_{\Omega} := \arg \min_{\beta \in \R^p}
\biggl \{ \|  Y- X \beta \|_n^2 + 2 \lambda \Omega ( \beta ) \biggr \} . $$
We will derive an oracle inequality for $\hat \beta$ for weakly decomposable norms $\Omega$,
 a notion introduced in Definition \ref{separable.definition}.

Recall that the $\ell_1$-norm is decomposable:
$\| \beta \|_1 = \| \beta_S \|_1 + \| \beta_{S^c} \|_1$ for all vectors $\beta$ and any set $S$. 
The triangle inequality of course holds for any norm $\Omega$ and so does the dual norm equality with
the uniform norm replaced by the dual norm 
$$\Omega_* (w) := \sup_{\Omega (\beta) \le 1 } | w^T \beta | . $$
We stress that the triangle inequality and dual norm equality fail to hold if we replace the norm by
powers of that norm. For example, the triangle inequality does not hold for
$\| \cdot \|_2^2$, which will mean the ridge regression penalty does not fall within
our framework. 
Returning to a general norm $\Omega$, it is not necessarily decomposable. Decomposability is however
very useful
for the derivation of oracle inequalities, an observation which was
discussed previously by 
\cite{vdG:2001}, 
\cite{SvdGICM} (where the property is called separability) and \cite{Wai11}.
Note that powers of norms can be decomposable, for example
$\| \beta \|_2^2 = \| \beta_S \|_2^2 + \| \beta_{S^c} \|_2^2$. However, the  required triangle inequality 
does not hold for $\| \cdot \|_2^2$. 

We will show now that decomposability is not a necessary condition for oracle results.
This was also realized by \cite{Bach:10}, although there the situation is restricted
to structured sparsity norms defined by sub-modular functions. 
We consider general norms $\Omega$, which are perhaps not decomposable, but
only {\it weakly decomposable} for certain index sets $S$, which means 
 that the norm $\Omega (\beta)$ of an arbitrary vector $\beta$ is always superior to the sum of
norms of $\beta_S$ and $\beta_{S^c}$. 

\begin{definition} \label{separable.definition} Fix some set $S$. We say that the norm
$\Omega$ is weakly decomposable if there exists a
norm $\Omega_{S^c}$  on $\R^{p-|S|} $ such that for all $\beta \in \R^p$,
$$\Omega (\beta) \ge \Omega (\beta_S ) + \Omega^{S^c} (\beta_{S^c} ) . $$
\end{definition}

\begin{definition} We say that $S$ is an allowed set (for $\Omega$) if
$\Omega$ is weakly decomposable for $S$. 
 \end{definition}

The best choice  for $\Omega^{S^c}$ is to take $\Omega^{S^c} ( \beta_{S^c} ) $ as large as possible
 (see also Section \ref{discussion}).
We identify $\beta_{S^c}$ with the $(p-|S|)$-vector $\{ \beta_j \}_{j \in S^c}$ and consider 
$\Omega^{S^c}$ as norm on $\R^{p-|S|}$ instead of $\R^p$. 
There may be no ``natural" extension to a norm on 
$\R^p$ (see Section \ref{Micchelli.example}  for an illustration), and
an extension is also not needed.

Observe that any norm is trivially (weakly) decomposable for the empty set and for the complete set 
$\{ 1 , \ldots , p \}$ containing the indices of the all variables.
Some examples, where we in particular discuss nontrivial choices of
 $S$, will be given in Section \ref{examples.section}.

We also extend the definition of $\ell_1$-eigenvalues and
effective sparsity to general weakly decomposable norms.

\begin{definition} Suppose $S$ is an allowed set.
 Let $L>0$ be some constant. The $\Omega$-eigenvalue (for $S$) is
 $$\delta_{\Omega} (L,S) := \min \biggl \{ \| X \beta_S - X{\beta_{S^c}} \|_n : \ 
 \Omega ( \beta_S) = 1 , \ \Omega^{S^c} ( \beta_{S^c} ) \le L  \biggr \}  .$$
 The $\Omega$-effective sparsity is 
 $$\Gamma_{\Omega}^2 (L,S) := { 1 \over \delta_{\Omega}^2 (L,S) } . $$
 \end{definition}
 
 The $\Omega$-eigenvalue $\delta_{\Omega} ( L ,S)$ depends on
 the choice of the norm $\Omega^{S^c}$, but we do not express this in our notation.
  It has a similar geometric interpretation as the $\ell_1$-eigenvalue:
 $\delta_{\Omega} ( L,S)$ is the distance between the sets
 $\{ X \beta_S: \ \Omega (\beta_S) = 1 \}$ and $\{ X \beta_{S^c}: \ \Omega^{S^c} (\beta_{S^c}) \le L \} $.
 The shape of these sets depends heavily on the norms $\Omega$ and $\Omega^{S^c}$.

 We will use the effective sparsity to bound the norm of $\beta_S$ in terms of $\| X \beta \|_n$, as 
 detailed in the following lemma. Here we use the ``cone condition" for $\Omega$.
 
 \begin{definition} Let $L >0$ be some constant, $S$ some allowed set and $\beta \in \R^p$ some vector.
 We say that $\beta$ satisfies the $(L,S)$-cone condition for $\Omega$ if
 $\Omega^{S^c} (\beta_{S^c}) \le L \Omega (\beta_S) $. 
 \end{definition}
 
 \begin{lemma} \label{rewrite.eigenvalue} 
 Suppose $S$ is an allowed set. 
 Then
 $$\delta_{\Omega} (L,S) =
 \min \biggl \{ { \| X \beta \|_n \over \Omega (\beta_S) } :\   \beta \ {\rm satisfies \ the }\ (L,S){\rm-cone \ condition }, \ \beta_S \not= 0   \biggr\}  $$
 and hence, for all $\beta$ that satisfy the $(L,S)$-cone condition,
 $$\Omega (\beta_S) \le \Gamma_{\Omega} (L,S) \| X \beta \|_n . $$ 
 \end{lemma}

 The ingredients for an oracle inequality are now:\\
 $\bullet$ the dual-norm equality,\\
 $\bullet$ the triangle inequality,\\
 $\bullet$ weak decomposability.\\
 In other words, the situation is as for the Lasso, but the
 decomposability property is weakened. 
  The dual norm of $\Omega $ is denoted by $\Omega_*$, 
  that is
  $$\Omega_* (w) := \sup_{\Omega(\beta) \le 1 } | w^T \beta | , \ w \in \R^p . $$
  We moreover let $\Omega_*^{S^c}$ be the dual norm of $\Omega^{S^c}$.

  \begin{theorem} \label{main.theorem} Let $\beta \in \R^p$ be arbitrary and let
  Let $S \supset \{ j : \beta_j \not= 0 \}  $ be an allowed set. Define
 $$\lambda^{S}  := \Omega_* \biggl ( (\epsilon^T X)_{S}/n  \biggr ) , \ \lambda^{S^c} :=
 \Omega_*^{S^c } \biggl ( (\epsilon^T X)_{S^c}/n \biggr ) . $$
 Suppose
 $$\lambda >  \lambda^{S^c} .$$
Define for some $0 \le \delta < 1$
 $$L_S:=  \biggl ( { \lambda + \lambda^{S} \over 
 \lambda - \lambda^{S^c} } \biggr ) \biggl ( { 1+ \delta \over 1- \delta} \biggr )  . $$
 Then
 $$\| X ( \hat \beta - \beta^0)  \|_n^2 + \delta (\lambda - \lambda^{S^c} ) \Omega^{S^c} (\hat \beta_{S^c})
 + \delta ( \lambda + \lambda^S) \Omega (\hat \beta_S - \beta)  $$ $$  \le \| X (\beta - \beta^0 )\|_n^2 +
\biggl  [ (1+ \delta ) (\lambda + \lambda^S ) \biggr ]^2 \Gamma_{\Omega}^2 (L_S , S) . $$
 
  \end{theorem}
  
  Theorem \ref{main.theorem} requires that $S \supset S_{\beta}$ is an allowed set.
   If, for values
   of $\beta$ that one considers as  good approximations of $\beta^0$,
    the smallest allowed set $S\supset S_{\beta}$ is much larger than $S_{\beta}$, then the penalty is
 simply not suited to describe the underlying sparsity structure.

  As a special case, one may take $\beta = \beta_0$ and $S_0 $ the smallest allowed set
  containing all non-zero $\beta_j^0 $ ($j=1 , \ldots , p $). However, the trade-off between approximation error
  $ \| X (\beta - \beta^0 )\|_n^2$ and 
  estimation error 
  $(\lambda + \lambda^S )^2 \Gamma_{\Omega}^2 (L_S , S)$
  will give better bounds.
  Theorem \ref{main.theorem} is sharp as the constant in front of the
  approximation error  $ \| X (\beta - \beta^0 )\|_n^2$ is one. 
  The choice $\delta =0$ is optimal if one only
  is interested in bounds for the prediction error $\| X ( \hat \beta - \beta^0 )\|_n^2 $.

 \section{Some examples}\label{examples.section}

  \subsection{The Lasso} \label{ell1.example} The $\ell_1$-norm $\Omega (\cdot) := \|\cdot \|_1$ is
(weakly) decomposable for all $S$, with $\Omega_{S^c}  = \Omega $, and $\Omega_* = \| \cdot \|_{\infty} $.
Hence, for all $\beta$ the set $S_{\beta}$ is an allowed set, that is, we can take $S =
S_{\beta}$ in Theorem \ref{main.theorem}. The choice $\delta =0$ then gives Theorem \ref{leastsquares}.
For $\delta >0$ however, we see that we also obtain a bound for 
$\| \hat \beta - \beta \|_1$, and hence for the $\ell_1$-estimation error error
$\| \hat \beta - \beta^0 \|_1 $. Here, one can use the triangle inequality
$\| \hat \beta - \beta^0 \|_1\|  \le \| \hat \beta - \beta \|_1+ \| \beta - \beta^0 \|_1 $, i.e., it again
involves a trade-off.

\subsection{Group Lasso} \label{group.example} 
Also the group Lasso norm $\| \cdot \|_{2,1}$ falls within the framework of decomposable norms.
Let $G_t \subset \{ 1 , \ldots , T\} $, 
$\cup_{t=1}^T  G_t = \{ 1 , \ldots , p \}$,
$G_1 \cap \cdots G_T = \emptyset$ be a partition of
$\{ 1 , \ldots , p \}$ into disjoint groups. The norm corresponding to the group Lasso penalty is
$$\Omega (\beta) := \| \beta \|_{2,1} :=
\sum_{t=1}^T \sqrt { | G_t | } \| \beta_{G_t} \|_2 , \ \beta \in \R^p . $$
It is (weakly) decomposable for $S= \cup_{t \in {\cal T}} G_t $ (${\cal T}$ being any subset
of $\{ 1, \ldots , T \}$), with $\Omega_{S^c} = \Omega$. 
Thus, we can take $S := \cup \{ G_t:\ \| \beta_{G_t} \|_2 \not= 0 \} $ as allowed set,
that is, as soon as $\beta_j \not= 0 $ for some $j \in G_t$, we take the whole group
of indexes $G_t$ into our allowed set $S$. 
The dual norm is
$$\Omega_* (w) := \| w \|_{2, \infty} := \max_{1 \le t \le T} \| w_{G_t} \|_2 / \sqrt {|G_t|} , \ w \in \R^p . $$
Let $X_{G_t} := \{ X_j \}_{j \in G_t} $ be the $n \times |G_t|$ design matrix of the variables in group $t$
($t=1 , \ldots , T$). Suppose that within groups the design is orthonormal, that is
$X_{G_t}^T X_{G_t} /n= I $ for all $t$. Then $\| X \beta_{G_t}\|_n = \| \beta_{G_t} \|_2 $ and
when $\epsilon \sim {\cal N} (0 ,  I)$, the random variables
$ \| (\epsilon^T X)_{G_t} \|_2^2  $ have a $\chi^2$-distribution with $|G_t|$ degrees of freedom.
Thus, 
$$\lambda_0^2 := \| (\epsilon^T X) \|_{2, \infty}^2 $$
is the maximum of $T$ normalized $\chi^2$-random variables. 
Invoking probability inequalities for such maxima, Theorem \ref{main.theorem}
then gives similar (but sharp) oracle results as those in \cite{lounici:11}
or \cite{BvdG2011}.

\subsection{General structured sparsity} \label{Micchelli.example}The following example describes a general structured sparsity
norm, as introduced by \cite{Micchelli:2010}. 
Let ${\cal A} \subset  [0, \infty )^p   $ be some convex cone, satisfying
${\cal A} \cup (0, \infty)^p  \not= \emptyset$, 
and
$$\Omega  (\beta ) := \Omega(\beta; {\cal A}) := \min_{a \in {\cal A}} {1 \over 2}
\sum_{j=1}^p \biggl ( { \beta_j^2 \over a_j } + a_j   \biggr ) . $$
Here we use the convention $0/0=0$. 
The assumption ${\cal A} \cup (0, \infty)^p  \not= \emptyset$ says that
there is an $a \in {\cal A}$ with all entries positive, so that for all $\beta$,
$\Omega (\beta )< \infty$.
It is shown in \cite{Micchelli:2010} that $\Omega$ is indeed a norm. 

Let ${\cal A}_S := \{ a_S : \ a \in {\cal A} \} $. 

\begin{definition} We call ${\cal A}_S$ an {\rm allowed}
set, if
$${\cal A}_S \subset {\cal A} . $$
\end{definition}
Thus we use the same terminology for sets in $\R^p$ (such as
${\cal A}_S$) and index sets $S$.

\begin{lemma}\label{allowed.lemma}
Suppose ${\cal A}_S$ is an allowed set. Then $S$ is allowed, that is if we
take
$$\Omega^{S^c}  (\beta_{S^c} ) = \Omega (\beta_{S^c} ; {\cal A}_{S^c}) , \ \beta_{S^c}  \in \R^{p-|S|} , $$
where ${\cal A}_{S^c} := \{ a_{S^c} : \ a \in {\cal A}  \} $, then the set $S$ is 
weakly decomposable for $\Omega$.
\end{lemma}

Note that ${\cal A}_{S^c}$ is a cone and that there always is an $a_{S^c} \in {\cal A}_{S^c}$
which has all entries positive except for those in ${\cal A}$.
Hence the restriction of $\Omega ( \cdot ; {\cal A}_{S^c} )$ to
$\{ \beta_{S^c} : \ \beta \in \R^p \}$ is a norm.
We do not require ${\cal A}_{S^c}$ to be an allowed set.

\begin{example} As in \cite{Micchelli:2010}, consider the convex cone
 $${\cal A}:= \{ a_1 \ge a_2 \ge \dots \ge a_p \ge 0 \} .$$
 The norm-penalty with norm $ \Omega ( \beta , {\cal A} ) $ then favors putting the
 last indexes equal to zero. Moreover, for any $s$, the set of the first $s$ indexes
 $\{ 1, \ldots , s \}$ is an allowed set.  A partition
 $\{ G_t \}_{t=1}^T $ is called {\it contiguous} if for all $t=1 ,\ldots , T-1$ and all
 $j \in G_t$ and $k \in G_{t+1}$ it holds that $j < k $.
 In \cite{Micchelli:2010} it is shown that for
 all $\beta$ there is a unique contiguous partition $\{ G_t \}_{t=1}^T$ of
 $\{1 , \ldots , p \}$ such that
 $$\Omega (\beta ; {\cal A} ) = \sum_{t=1}^T \sqrt {|G_t|} \| \beta_{G_t} \|_2 . $$
 \end{example}

 We now return to the general norm $\Omega( \cdot ; {\cal A})$. Its dual norm is
 $$\Omega_* ( w ; {\cal A} ) = \max_{a \in {\cal A} (1) } 
 \sqrt { \sum_{j=1}^p a_j w_j^2 } , \ w \in \R^p,$$
 where ${\cal A} (1) := \{ a \in {\cal A}:\ \| a \|_1 = 1 \}$. 
 A similar expression holds for the dual norm $\Omega_*^{S^c}$ of
 $\Omega^{S^c}$.
 
 \cite{PontilMaurer}
 provide moment inequalities for $\Omega_* ( \epsilon^T X; {\cal A})$.
 They show that when $\epsilon \sim {\cal N} (0, I)$, then
$$\EE  \Omega_* ( \epsilon^T X ; {\cal A} ) /n \le \lambda_{\epsilon}
,$$
with
$$\lambda_{\epsilon} := \sqrt { 8 \over n  }  \biggl (2+ \sqrt { \log | {\rm extreme \ points \ of\ } {\cal A} (1) |  }\biggr ) 
\sqrt { \sum_{i=1}^n \Omega_*^2 ( x_i  ; {\cal A} ) \over n}  , $$
where $x_i = (x_{i,1} , \ldots , x_{i,p} )$ is the $i$-th row of $X$. 
Using concentration of measure (\cite{Talagrand:95}) , this can be turned into a suitable
probability inequality. Again, the results can be applied
to $\Omega_*^{S^c}$ as well.

 The $\ell_1$-norm is a special case of the structured sparsity norm,
 with ${\cal A} = [0, \infty)^p $. 
 
 The norm $\| \cdot \|_{2,1} $ corresponding to the group Lasso,
 as described in Subsection \ref{group.example} is also a special case, with
 $${\cal A} := \{ a \in [0, \infty)^p   \ {\rm is \ constant \ within  \ groups }  \}. $$

 \subsection{A trivial example}
\label{trivial.example} A trivial example is the norm
$$\Omega_G(\beta) := \sqrt {|G|} \| \beta_G\|_2 + \| \beta_{G^c} \|_1 , $$
which is a special case of the group Lasso norm, with $n-|G|+1$ groups,
namely, the group $G$ and $n- |G| $ groups $\{ j \}_{j\notin G} $, each containing only one
element. It is weakly decomposable for each $S\supset G$ with $\Omega^{S^c} = \| \cdot \|_1 $.
We will invoke this example mainly for facilitating our discussion of
the relation between $\Omega$-eigenvalues (see Section \ref{eigenvalues.section}).
 
 \subsection{Overlapping groups}  \label{overlapping.example}
 In this example, we consider a norm corresponding to the group Lasso
 with overlapping groups (\cite{Jacob:2009}). Let $\{G_t\}_{t=1}^T$ be subsets of $\{ 1 , \ldots , p \}$,
 with $\cup_{t=1}^T G_t = \{ 1 , \ldots , p \}$, and define
 $$\Omega_{\rm overlap} (\beta) := \min\biggl \{ \sum_{j=1}^T \| b_t \|_2: \ (b_t)_{G_t^c} =0  \ \forall \ t , \ \sum_{t=1}^T b_t = \beta \biggr \} . $$
 The paper \cite{Jacob:2009} shows that $\Omega_{\rm overlap}$  is indeed a norm.
 However, as such $\Omega_{\rm overlap}$ is not weakly decomposable for useful candidate sets $S$. 
 On the other hand by a reparametrization with parameters $\{ b_t \}_{t \in T} $, we can reformulate the overlapping group
 Lasso problem into a group Lasso problem with non-overlapping groups. To see this, note that
 $$ {X}\beta= \sum_{t=1}^T { X} b_t , \ \sum_{t=1}^T b_t = \beta . $$
 Thus, the overlapping group Lasso estimator is 
 $\hat \beta = \sum_{t=1}^T \hat b_t $, 
 where
 $$\{ \hat b_t \}_{t=1}^T := \arg \min_{ \{ b_t \}_{t=1}^T : \ ({b_t})_{G_t^c} = 0  \ \forall \ t } \biggl \{ \| Y - \sum_{t=1}^T {X} b_t \|_n^2 +
 2 \lambda \sum_{t=1}^T \| b_t \|_2  \biggr \} . $$
 The augmented model has $\tilde p := \sum_{t=1}^T |G_t |$ parameters $\{ b_{j,t} :\ j \in G_t \}_{t=1}^T $ and the augmented groups
 are $\tilde G_t := \{ (j,t):\ j \in G_t \}$ ($t=1, \ldots , T$), which are by definition non-overlapping. 
 However, in the augmented design matrix
 $$\tilde {X} := \biggl\{ \{ {X}_j : \ j \in G_t \} \biggr \} _{t=1}^T $$
the column ${X}_j$ appears $N_j := \sum_{t=1}^T {\rm l } \{ j \in G_t \}$ times ($j=1 , \ldots , p$). 
 Although such repetitions are not a problem for the Lasso (see Remark \ref{repeat.remark}), the implications for the group Lasso are not so clear.

%

\section{Comparing $\Omega$-eigenvalues}\label{eigenvalues.section}

The question arises to what extend using a norm-penalty with norm
$\Omega$ different from the
$\ell_1$-norm results in better oracle inequalities. This partly depends
on the behavior of the dual norm, a topic we briefly discuss in Section \ref{discussion}. It also depends on the
behavior of the $\Omega$-eigenvalues, which is the theme of the present section.

Fix a set $S$ and
 consider again the norm  $\Omega_S$-defined in Section
 \ref{trivial.example}:
 $$\Omega_S (\beta) = \sqrt {|S|} \| \beta_S \|_2 + \| \beta_{S^c} \|_1 . $$
 This norm is decomposable for $S$ with $\Omega^{S^c} = \| \cdot \|_1 $.
 The $\Omega_S$-eigenvalue
 $\delta_{\Omega_S} (L,S) $ is the distance between the
 contour of the ellipse $ \{ X \beta_S : \ \| \beta_S \|_2 = 1 /\sqrt {|S|}\}$ and the 
 convex hull including interior $\{  X \beta_{S^c}: \ \| \beta_{S^c} \|_1 \le L \} $. 
 
 \begin{remark}
 In fact, $\delta_{\Omega_S } ( L, S)$ is in part easy to compute: for fixed
 $\beta_{S^c}$ one calculates
 $$\min_{\| \beta_S \|_2^2 = 1/{|S|} } \| X \beta_S - X \beta_{S^c} \|_n^2 := {\cal R}^2 (\beta_{S^c}) . $$
 This is a quadratic minimization problem with quadratic restriction, which can be 
 solved using Lagrange calculus. The more difficult part is to find the minimizer of 
 $ {\cal R}^2 (\beta_{S^c}) $ over all $\| \beta_{S^c} \|_1 \le L $. 
 \end{remark}
 
 In \cite{BvdG2011}, $|S| \times \delta_{\Omega_S}^2 (L,S) $ is called the
 {\it adaptive restricted eigenvalue} (because it occurred there in conjunction with the adaptive Lasso).
 
 Recall that $\delta(L,S)$ is the $\ell_1$-eigenvalue.
 Since $\| \beta_S \|_1 \le \sqrt {|S|} \| \beta_S \|_2$, one easily checks that
 $$\delta (L,S) \ge \delta_{\Omega_S} (L,S) , $$
 i.e., the $\ell_1$-eigenvalue $\delta(L,S)$ is better behaved than the
 $\Omega_S $-eigenvalue $\delta_{\Omega_S} (L,S)$. 
 
 Consider now the structured sparsity norm
 $\Omega (\cdot ; {\cal A})$ introduced in Section \ref{Micchelli.example}.
 By Lemma \ref{allowed.lemma}, we know that under the condition that ${\cal A}_S $ is allowed, the norm
 $\Omega (\cdot ; {\cal A})$ is weakly decomposable for $S$ with
 $\Omega^{S^c} (\beta_{S^c} ) = \Omega ( \beta_{S^c} ;  {\cal A}_{S^c} ) $.

We note that
$$\Omega(\beta; {\cal A} ) \ge \| \beta \|_1 , $$
and $\iota_S \in {\cal A}$, where $\iota$ is the constant vector $\iota:= (1 , \cdots, 1)$,
then
$$\Omega (\beta_S;{\cal A} ) \le \sqrt {|S|} \| \beta_S \|_2 . $$
In other words, $\Omega ( \cdot , {\cal A})$ intermediates the $\ell_1$- and
$\ell_2$-norm. 

\begin{lemma} \label{compare.lemma}
Suppose ${\cal A}_S$ is allowed and that 
$\iota_S \in {\cal A}$, where $\iota$ is the constant vector $\iota:= (1 , \cdots, 1)$.
Then for all $L >0$, 
$$\delta_{\Omega} (L, S) \ge \delta_{\Omega_S} ( L, S) . $$
\end{lemma}

It follows that 
$$\Gamma^2(L,S )\le \Gamma_{\Omega_S}^2 (L,S)  . $$
and more generally, under the conditions of
Lemma \ref{compare.lemma} 
$$\Gamma_{\Omega}^2  (L, S) \le \Gamma_{\Omega_S}^2 (L,S).$$
The $\Omega$-effective sparsity $\Gamma_{\Omega}^2 (L,S) $  is in general not comparable to
the $\| \cdot \|_1$-effective sparsity $\Gamma^2 (L,S)$ for the $\ell_1$-norm $\| \cdot \|$. This is only partly due to the fact that the
cone condition for $\Omega$  and the cone condition for $\| \cdot \|_1$ are not comparable.
We finally note that the restricted eigenvalue (see \cite{bickel2009sal}) is in between 
$|S| \delta_{\Omega_S}^2 (L,S)$ and $|S| \delta^2 (L,S)$, and that the $\Omega$-eigenvalue $\delta_{\Omega} (L,S)$ is not comparable to
the restricted eigenvalue either, which is now solely due to the incomparability of the cone conditions. 

\section{Discussion} \label{discussion}

We have shown that sparsity oracle properties hold for the least squares
estimator with separable norm-penalty.  There are a few issues that can be addressed here.

First of all, the choice of a norm other than $\| \cdot \|_1$ can be inspired by the practical use: the estimator
may have a better interpretation. On the other hand, it may be harder to compute.

The second point is that with another norm,
the dual norm may better behaved than with the $\ell_1$norm. This is the case for for instance the group Lasso,
which wins in certain cases  from the Lasso by a $\log p$-term. In this paper,
we have not discussed in detail the properties of the dual 
$\Omega_* \left ( (\epsilon^T X )_S \right) $ or $\Omega_*^{S^c}\left  ((\epsilon^T X)_{S^c} \right ) $
to avoid digressions. General results can be found in \cite{PontilMaurer}. 
Larger norms have smaller dual norms, that is if $\Omega (\beta)
\ge \tilde \Omega (\beta)$ for all $\beta$, then
$\Omega_* (w) \le \tilde \Omega_* (w) $ for all $w$. Note that Theorem
\ref{main.theorem} gives bounds for the $\Omega$-error of $\hat \beta_{S_0}$, so not only its
dual norm is smaller than that of $\tilde \Omega$, but also the bound holds for the $\tilde \Omega$-error.
In particular, this comparison can be made between  the structured sparsity norm 
$\Omega (\cdot ; {\cal A}) $ defined in Section \ref{Micchelli.example} and
the $\ell_1$-norm, because
$\Omega (\beta ; {\cal A} ) \ge  \| \beta \|_1 $ for all $\beta$.
Note further that Theorem \ref{main.theorem} also 
involves $\Omega^{S^c}$ and
its dual $\Omega_*^{S^c} $, and that its result can be optimized
by taking the largest possible choice for $\Omega^{S^c}$ (which will then also
optimize the $\Omega$-eigenvalue).

Of course, the prize to pay for using a norm different from
$\ell_1$ is that it may only be weakly decomposable for relatively large sets $S$. 
That is, one should choose a norm that corresponds to a priori knowledge on the
sparsity structure.

It is to be noted further that with invoking the dual norm equality one might not exploit in full
the structure of the problem. More refined techniques are given in for example
\cite{Lederer:11}. 

In cases where the penalty involves a ``smoothness" norm
(for example a Sobolev norm), the philosophy is again
different. In the classical setup, such a penalty is invoked for
establishing (non-adaptive) smoothness only. In more recent settings,
the aim is to obtain both sparsity and smoothness. An example,
concerning the high-dimensional additive model,
is in \cite{meier08addmodeling}. There, the issue of 
decomposability,
comes up as well. Oracle results are derived using a penalty
that is not only sparsity decomposable  but also ``smoothness"
decomposable (see also \cite{BvdG2011}, Section 8.4.5). 

Finally, the oracle results can be extended to loss functions other than least squares
(for example in the spirit of \cite{vandeG08} or \cite{Wai11}). Sharp oracle
results are discussed in \cite{vdG:2013}.
For the quasi-likelihood loss with canonical link function,
the dual-norm argument can again be used.
For other cases  this
argument generally has to be replaced. Here, tools from empirical
process theory can be invoked (such as those outlined in
\cite{BvdG2011}, Chapter 8).

\section{Proofs} \label{proofs.section}

 {\bf Proof of Lemma \ref{rewrite.eigenvalue}.} Let $\beta \in {\cal C}:= \{ \Omega^{S^c} ( \beta_{S^c}  ) \le L \Omega ( \beta_S) \not= 0\} $.
Write
 $$\tilde \beta_S := { \beta_S \over \Omega (\beta_S ) } ,\ 
 \tilde \beta_{S_c} := {\beta_{S^c} \over \Omega (\beta_S ) } . $$
 Then $\Omega ( \tilde \beta_S) = 1$ and $\Omega^{S^c} ( \tilde
 \beta_{S^c} ) \le L $, and hence
 $${ \| X \beta \|_n  \over \Omega(\beta_S) } =
 \| X \tilde \beta_S + X \tilde \beta_{S^c} \|_n .$$
 It follows that
 $$\min_{\beta \in {\cal C} } { \| X \beta \|_n  \over \Omega(\beta_S) } = \delta_{\Omega}  (L,S). $$
 \hfill $\sqcup \mkern -12mu \sqcap$
 
 The next lemma shows why convexity of the penalty is important.
 The result can be extended to loss functions other than quadratic loss, see the rejoinder in the discussion paper \cite{vdG:2013}. 
 
  \begin{lemma} \label{convexpenalty.lemma} Let ${\cal B}$ be a convex subset of $\R^p$ and ${\rm pen}: {\cal B}  \rightarrow \R$ be a convex penalty.
  Let moreover
 $$ \hat \beta := \arg \min_{\beta \in {\cal B}} \| Y - X \beta \|_n^2 + 2{\rm pen} (\beta) . $$
 Then for every $\beta \in {\cal B}$
 $$(Y- X \hat \beta )^T X(\beta - \hat \beta )/n  + {\rm pen} (\hat \beta) \le {\rm pen} (\beta) . $$
  \end{lemma}
  
  {\bf Proof.} Fix $\beta \in {\cal B}$ and define for $0 < \alpha \le 1$, 
  $$ \hat \beta_{\alpha} := (1- \alpha) \hat \beta + \alpha \beta . $$
  
 We have
  $$\| Y - X \hat \beta \|_n^2 + 2{\rm pen} (\hat \beta) \le \| Y - X \hat \beta_{\alpha} \|_n^2 + 2{\rm pen} ( \hat \beta_{\alpha} )  $$
  $$ \le \| Y - X \hat \beta_{\alpha} \|_n^2 + 2(1- \alpha )  {\rm pen} ( \hat \beta) + 2\alpha {\rm pen} (\beta ) $$
  where we used the convexity of the penalty. 
  It follows that
  $${ \| Y - X \hat \beta \|_n^2 - \| Y - X \hat \beta_{\alpha} \|_n^2 \over \alpha } + 2{\rm pen} (\hat \beta)\le 2{\rm pen} (\beta) . $$
 But clearly
  $$ \lim_{\alpha \downarrow 0 } { \| Y - X \hat \beta  \|_n^2 - \| Y- X \hat \beta_{\alpha} \|_n^2 \over \alpha } =
  2 ( Y - X \hat \beta )^T X(\beta - \hat \beta ) /n . $$
  
  \hfill $ \sqcup \mkern -12mu \sqcap$
  
  {\bf Proof \ref{main.theorem}.} Let us write for $v,w \in \R^n$,
  $$(v,w) := v^T w / n . $$
  Fix some $\beta \in \R^p$  and let $S \supset  \{ j: \ \beta_j \not= 0 \} $ be an allowed set
  If 
  $$( X( \hat \beta - \beta^0) , X ( \hat \beta -\beta ) )_n \le -( \delta (\lambda +  \lambda^S ) \Omega ( \hat \beta_S - \beta ) + \delta (\lambda - \lambda^{S^c} )  \Omega (\hat \beta_{S^c}))$$
  we find
  $$\| X ( \hat \beta - \beta^0 ) \|_n^2 + \delta (\lambda +  \lambda^S ) \Omega ( \hat \beta_S - \beta ) + \delta (\lambda - \lambda^{S^c} )  \Omega (\hat \beta_{S^c})
- \| X (\beta - \beta^0 ) \|_n^2 $$
  $$ = \delta (\lambda +  \lambda^S ) \Omega ( \hat \beta_S - \beta ) + \delta (\lambda - \lambda^{S^c} )  \Omega (\hat \beta_{S^c})
   - \| X ( \beta - \hat \beta ) \|_n^2 + 2 ( X( \hat \beta - \beta^0) , X ( \hat \beta - \beta ) )_n \le 0 . $$
  Hence, then we are done.
  
  Suppose now that
  $$( X( \hat \beta - \beta^0) , X ( \hat \beta - \beta ) )_n \ge -  (\delta (\lambda +  \lambda^S ) \Omega ( \hat \beta_S - \beta ) + \delta (\lambda - \lambda^{S^c} )  \Omega (\hat \beta_{S^c}))  .$$
  By Lemma \ref{convexpenalty.lemma} we have
  $$((Y- X \hat \beta ), X(\beta - \hat \beta ))_n  + \lambda \Omega (\hat \beta) \le \lambda \Omega (\beta) , $$
  or
  $$(X (\hat \beta -\beta^0 ), X(\hat \beta - \beta ))_n  + \lambda \Omega (\hat \beta) \le (\epsilon, X(\hat \beta - \beta))_n + \lambda \Omega (\beta) . $$
  By definition of the dual norm,
  $$(\epsilon, X(\hat \beta - \beta))_n = (\epsilon, X (\hat \beta_S - \beta))_n +
  (\epsilon, X (\hat \beta_{S^c} - \beta))_n \le \lambda^S \Omega ( \hat \beta_S - \beta ) +
  \lambda^{S^c} \Omega^{S^c} ( \hat \beta_{S^c} ) . $$
  Thus
  $$ (X (\hat \beta -\beta^0 ), X(\hat \beta - \beta ))_n  + \lambda \Omega (\hat \beta) \le 
   \lambda^S \Omega ( \hat \beta_S - \beta ) +
  \lambda^{S^c} \Omega^{S^c} ( \hat \beta_{S^c} ) + \lambda \Omega (\beta) . $$
  By the weak decomposability of $\Omega$ and the triangle inequality, this implies
  \begin{equation}\label{tea.equality}
  (X (\hat \beta -\beta^0 ), X(\hat \beta - \beta ))_n  + (\lambda - \lambda^{S^c} )  \Omega^{S^c} (\hat \beta_{S^c}) \le 
   (\lambda +  \lambda^S ) \Omega ( \hat \beta_S - \beta )  . 
   \end{equation} 
 Since $  (X (\hat \beta -\beta^0 ), X(\hat \beta - \beta ))_n \ge - (\delta (\lambda +  \lambda^S ) \Omega ( \hat \beta_S - \beta ) + \delta (\lambda - \lambda^{S^c} )  \Omega^{S^c} (\hat \beta_{S^c})) $ this gives
 $$\Omega^{S^c} ( \hat \beta_{S^c} ) \le L_S \Omega (\hat \beta_S - \beta) . $$
 We now insert Lemma \ref{rewrite.eigenvalue}, which gives
 \begin{equation}\label{coffee.equation}
 \Omega(\hat \beta_S - \beta ) \le \Gamma_{\Omega} ( L,S) \| X ( \hat \beta - \beta ) \|_n 
 \end{equation}
 and continue with inequality (\ref{tea.equality}):
 $$ (X (\hat \beta -\beta^0 ), X(\hat \beta - \beta ))_n  + (\lambda - \lambda^{S^c} )  \Omega (\hat \beta_{S^c}) + \delta ( \lambda + \lambda^S) 
 \Omega (\hat \beta_S - \beta) $$ $$ \le 
   [ (1+\delta) (\lambda + \lambda^S) ]\Gamma_{\Omega} (L_S , S)   \| X ( \hat \beta - \beta ) \|_n  $$
   $$ \le {1 \over 2} \biggl [  (1+ \delta ) (\lambda + \lambda^S  ) \biggr ]^2 \Gamma_{\Omega}^2  (L_S , S)  + {1 \over 2} \| X ( \hat \beta - \beta ) \|_n^2 . $$
   Since
   $$
   2 ( X( \hat \beta - \beta^0) , X ( \hat \beta - \beta ) )_n = \| X ( \hat \beta - \beta^0 ) \|_n^2  - \| X (\beta - \beta^0 ) \|_n^2 +
   \| X ( \beta - \hat \beta ) \|_n^2  ,$$
   we obtain
   $$\| X (\hat \beta - \beta^0 ) \|_n^2 + 2(\lambda - \lambda^{S^c} )  \Omega (\hat \beta_{S^c}) +  2 \delta ( \lambda + \lambda^S)   \Omega( \hat \beta_S - \beta) $$ $$
\le \| X (\beta - \beta^0 ) \|_n^2 + \biggr [ (1+ \delta)  (\lambda + \lambda^S  )\biggr ] ^2 \Gamma_{\Omega}^2  (L_S , S) . $$
  \hfill $\sqcup \mkern -12mu \sqcap$
 
 {\bf Proof of Lemma \ref{allowed.lemma}.} 
Note that for any $a$ and $\beta$
$$ {1 \over 2}
\sum_{j=1}^p \biggl ( { \beta_j^2 \over a_j } + a_j   \biggr ) =
{1 \over 2}
\sum_{j\in S} \biggl ( { \beta_j^2 \over a_j } + a_j   \biggr )+
{1 \over 2}
\sum_{j\in S^c } \biggl ( { \beta_j^2 \over a_j } + a_j   \biggr ). $$
Hence,
writing
$$a (\beta) := \argmin_{a \in {\cal A}} {1 \over 2}
\sum_{j=1}^p \biggl ( { \beta_j^2 \over a_j } + a_j   \biggr ) , $$
we have
$$\Omega (\beta) := {1 \over 2}
\sum_{j=1}^p \biggl ( { \beta_j^2  \over a_j (\beta)} +a_j (\beta)  \biggr ) $$
 $$ =
{1 \over 2}
\sum_{j\in S} \biggl ( { \beta_j^2 \over a_j (\beta) } + a_j  (\beta)  \biggr )+
{1 \over 2}
\sum_{j\in S^c }^p \biggl ( { \beta_j^2 \over a_j (\beta) } + a_j  (\beta)  \biggr )$$
$${1 \over 2}
\sum_{j=1}^p  \biggl ( { \beta_{j,S}^2 \over a_{j,S} (\beta) } + a_{j,S}  (\beta)  \biggr )+
{1 \over 2}
\sum_{j=1 }^p \biggl ( { \beta_{j, S^c}^2 \over a_{j, S^c} (\beta) } +
 a_{j, S^c  }  (\beta)  \biggr )$$
 $$\ge \min_{a_S \in {\cal A}_S} {1 \over 2}
\sum_{j=1}^p  \biggl ( { \beta_{j,S}^2 \over a_{j, S} } + a_{j, S}    \biggr )+
\min_{a_{S^c} \in {\cal A}_{S^c} } {1 \over 2}
\sum_{j=1 }^p \biggl ( { \beta_{j, S^c}^2 \over a_{j, S^c}  } +
 a_{j, S^c  }    \biggr )$$
 $$\ge \Omega ( \beta_S  ) + \Omega^{S^c}  ( \beta_{S^c}) .$$
 \hfill $\sqcup \mkern -12mu \sqcap$
 
 {\bf Proof of Lemma \ref{compare.lemma}.} Suppose $\beta$ satisfies the $(L,S)$-cone condition for $\Omega$:
$$\Omega^{S^c}  (\beta_{S^c} ) \le L \Omega ( \beta_S), $$
then also
$$\Omega_S^{S^c} ( \beta_{S^c}) =  \| \beta_{S^c} \|_1 \le \Omega^{S^c}  (\beta_{S^c} ) \le L
\Omega (\beta_S) \le L \sqrt {|S|} \| \beta_S \|_2 = L \Omega_S (\beta_S)  ,$$
where in the last inequality we used $\iota_S \in {\cal A}$.
Hence, $\beta$ satisfies the $(L,S)$-cone condition for $\Omega_S$.
But then
$$\Omega (\beta_S) \le \Omega_S (\beta_S) \le 
{ \| X \beta \|_n \over \delta_{\Omega_S}  ( L,S)}. $$
The result now follows from Lemma \ref{rewrite.eigenvalue}.

\hfill $\sqcup \mkern -12mu \sqcap$

\bibliographystyle{plainnat}
\bibliography{reference}

\end{document}